\documentclass[11pt]{article}
\usepackage[dvips]{graphics,graphicx}
\DeclareGraphicsExtensions{.eps}

\usepackage{a4}

\def\ds{\displaystyle}

\def\id{\mbox{\small 1}\!\!1}

\def\nn{\nonumber}
\def\re{\mbox{\sf Re}}

\def\stack{\stackrel}

\def\qed{\hbox{${\vcenter{\vbox{                        
   \hrule height 0.4pt\hbox{\vrule width 0.4pt height 6pt
   \kern5pt\vrule width 0.4pt}\hrule height 0.4pt}}}$}}

\def\R{\bf\sf R}
\def\Z{\bf\sf Z}
\def\C{\bf\sf C}

\newtheorem{theorem}{Theorem}
\newtheorem{lemma}{Lemma}

\newtheorem{definition}{Definition}

\begin{document}

\begin{center}
  
{\bf \Large A New Proof of a Theorem in Analysis by Generating Integrals 
and Fractional Calculus}

\bigskip

{\large S.C. Woon}

\medskip

Trinity College, University of Cambridge, Cambridge CB2 1TQ, UK

s.c.woon@damtp.cam.ac.uk

MSC-class Primary 28A25, 11M06; Secondary 26A33

Keywords: Theories of Integration; the Riemann zeta function\\
Analytical Number Theory; Fractional Calculus

December 23, 1998

\bigskip

{\bf \large Abstract}

\end{center}

The idea of generating integrals analogous to generating functions
is first introduced in this paper.  A new proof of the well-known
Finite Harmonic Series Theorem in Analysis and Analytical Number Theory
is then obtained by the method of Generating Integrals and Fractional
Calculus.  A generalization of the Riemann zeta function up to
non-integer order is derived.

\bigskip

\section{The Finite Harmonic Series Theorem}

\begin{definition}
The finite harmonic series is
$$ h(n) = \sum_{k=1}^{n}\frac{1}{k} = 1 + \frac{1}{2} + \cdots +
\frac{1}{n}\;.$$
\end{definition}

\begin{theorem}[{\bf The Finite Harmonic Series Theorem}]\quad\\
It is well known \cite[p. 16, (1.7.9)]{bateman} that 
\begin{eqnarray}
h(n) &=& \psi(1\!+\!n)+\gamma \label{e:harmthm1}\\
&=& \log n + \gamma + O(1/n) \label{e:harmthm2}
\end{eqnarray}
where the digamma function $\ds \psi(z) = \frac{d}{dz}
\log\left(\Gamma(z)\right) = \frac{\Gamma'(z)}{\Gamma(z)}$, the Euler
constant $\gamma = -\,\Gamma'(1) = 0.577215\cdots$, and the prime
denotes differentiation.
\end{theorem}

The first part (\ref{e:harmthm1}) of the Theorem can be proved by
differentiating the recurrence relation of
$\Gamma(n\!+\!1)=n\,\Gamma(n)$ and summing over $n$
\cite[p. 256]{abramowitz}. The second part (\ref{e:harmthm2}) can be
proved by using the Euler-Maclaurin summation formula
\cite[p. 629]{proof}.

In this paper, we shall present an alternative elementary proof
\cite{woon_harmonic} of the first part (\ref{e:harmthm1}) of the
Theorem by using {\bf Fractional Calculus} and
the idea of {\bf generating integrals}.

\section{Introduction \& Motivation of Fractional Calculus}

\begin{quote}
{\it ``Thus it follows that $d^{1/2}x$ will be equal to \dots from which
one day useful consequences will be drawn.''} --- Leibniz in a letter
\cite{leibniz} to L'Hospital.
\end{quote}

\begin{quote}
{\it ``When $n$ is an integer, the ratio $d^n p$, $p$ a function of $x$, to
$dx^n$ can always be expressed algebraically. Now it is asked: what
kind of ratio can be made if $n$ is a fraction?''} --- Euler \cite{euler}.
\end{quote}

\begin{quote}
{\it ``The idea of an integral or derivative, of arbitrary non-integral
  order, was introduced into analysis by Liouville and Riemann.  Such
  integrals and derivatives may be, and have been by different
  writers, defined in a variety of manners, and different systems of
  definitions may be the most useful in different fields of
  analysis.''} --- Hardy and Littlewood \cite{hardylwood1}.
\end{quote}

In this paper, we shall consider $D^\sigma, \;\sigma\in\R$, an
operator of non-integer order, a notion first pondered upon by Leibniz
\cite[pp. 301-302]{leibniz}, Euler \cite[p. 55]{euler}, Lagrange
\cite{lagrange}, Laplace \cite[p. 85 and p. 186]{laplace}, Fourier
\cite{fourier} and Abel \cite{abel} during the late 17th century to
early 19th century.  We shall briefly review several known and equally
valid definitions for $D^\sigma$, and then focus on one of the
definitions, the Riemann-Liouville (R-L) Fractional Calculus.  The
foundation of R-L Fractional Calculus was laid by Riemann
\cite{riemann} and Liouville \cite{liouville} in the late 19th
century, and then subsequently developed by Cayley \cite{cayley},
Laurent \cite{laurent}, Heaviside \cite{heaviside}, Hardy and
Littlewood \cite{hardylwood1, hardy1, hardylwood2, hardylwood3}, and
many others.  It was largely regarded as a mathematical curiosity
until only recently when Mandelbrot, the discoverer of Fractals, found
an application of the R-L Fractional Calculus in the Brownian motion
in a fractal medium, and speculated a possible connection between the
analysis of Fractional Calculus and the geometry of Fractals
\cite{mandelbrot}.

We shall also introduce the idea of generating integrals by
analogy to generating functions.  As we shall see, just
as a certain generating function is useful for generating a certain
desired sequence of numbers, a certain generating integral is
similarly useful. While a generating function $f(z)$ generates a
sequence of numbers $\{p_n\}$ in the coefficients of the terms of
different orders in its power series expansion $\sum_n p_n z^n$, a
generating integral generates a sequence of numbers $q_n$ in the
coefficient of a term in the result of an $n$-fold integration.

However, a generating integral has one unique advantage over a
generating function. The R-L Fractional Calculus can be used to
analytically extend a generating integral of iteration order
$n\in\Z^+$ to order $\rho\in\R$. The result is that the sequence of
numbers $\{q_n\}$ is in turn analytically extended to a function
$q(\rho), \;\rho\in\R$.

We shall then show how the Riemann-Liouville Fractional Calculus and
the idea of generating integrals can be used to prove the well-known
Finite Harmonic Series Theorem.

\section{Differential-Integral Operator $D^n$}

\begin{definition}
Let the operator $D^n_{\!x|a}, \;n\in\Z$, acting on a function $f$ at
the point $x$ be defined as
\begin{equation}
D^n_{\!x|a} f(x) = \left\{\begin{array}{ll}
\ds \frac{d^n}{dx^n} f(x) & (n>0)\\
\\
\ds\Big. f(x) & (n=0)\\
\\
\ds\Bigg. \int_a^x f(\hat{x}) (d\hat{x})^{-n} & (n<0)
\end{array}\right.
\label{e:Dndef}
\end{equation}
where the n-fold integration is defined inductively as
$$\int_a^x f(\hat{x}) (d\hat{x})^{-n} = 
\underbrace{\int^x_a\!\int^{x_n}_a\!\!
  \int^{x_{n-1}}_a\!\!\!\!\!\!\!\!\cdots\!\int^{x_3}_a
  \!\!\int^{x_2}_a\!}_{n \mbox{\footnotesize-times}} f(x_1) \;
dx_1\,dx_2 \cdots dx_{n-1}\,dx_n\;.$$
\end{definition}

As an example, consider $f(x) = x^m$.
\begin{equation}
D^n_{\!x|a} x^m = \left\{\begin{array}{ll}
0 & (m\ge 0, \;m<n)\\
\\
\ds\Bigg. \frac{m!}{(m\!-\!n)!}\,x^{m-n} =
\frac{\Gamma(1\!+\!m)}{\Gamma(1\!+\!m\!-\!n)}\,x^{m-n} & (m\ge 0,
\;m\ge n)\\
\\
\ds\Bigg. (-1)^n \frac{(|m|\!-\!1\!+\!n)!}{(|m|\!-\!1)!}
\,\hat{x}^{m-n} \Big|_a^x 
& (m\le 0, \;m<n)\\
\ds\Bigg. \quad = \lim_{\epsilon\to 0}
\frac{\Gamma(1\!+\!m\!+\!\epsilon)}{\Gamma(1\!+\!m\!+\!\epsilon\!-\!n)}\,
\hat{x}^{m-n}\Big|_a^x\\
\\
\ds\Bigg. \int_a^x \left( \int_a^{\hat{x}} \tilde{x}^m (d\tilde{x})^{|m|}
\right) (d\hat{x})^{(m-n)} & (m\le 0, \;m\ge n)\\
\ds\Bigg. \quad = (-1)^{m+1} \frac{1}{(|m|\!-\!1)!} \int_a^x \log \hat{x}
(d\hat{x})^{(m-n)}
\end{array}\right.
\label{e:Dnxdef}
\end{equation}
where $\ds \hat{x}^{m-n}\Big|^x_a = (x^{m-n} - a^{m-n})$.

\smallskip

If we tabulate $D^n_{\!x|a} x^m$ for $n,m\in\Z$ and omit the constant
terms containing $a$, we can observe a pattern emerges as in Table 1.

\begin{table}[hbt]
\begin{center}
$\begin{array}{|r||r|r|r|r|r|r|r|}
\hline
  \multicolumn{8}{|c|}{D^n_{\!x} x^m} \\
  \hline
  \mbox{\small\it m}\!\backslash\!\mbox{\small\it n} & -3 & -2 & -1 & 0
  & 1 & 2 & 3 \\
  \hline \hline 2 & \mbox{\small 2!/5!}\;x^5 & \mbox{\small
    2!/4!}\;x^4 & \mbox{\small 2!/3!}\;x^3 & x^2 & \mbox{\small 2!}\;x
  & \mbox{\small 2!} & \mbox{\small 0}  \\
  \cline{7-7} 1 & \mbox{\small 1/4!}\;x^4 & \mbox{\small 1/3!}\;x^3 &
  \mbox{\small 1/2!}\;x^2 & x & \mbox{\small 1} & \mbox{\small 0} &
  \mbox{\small 0} \\
  \cline{6-6} 0 & \mbox{\small 1/3!}\;x^3 & \mbox{\small 1/2!}\;x^2 &
  x & \mbox{\small 1} & \mbox{\small 0} & \mbox{\small 0} & \mbox{\small 0} \\
  \cline{2-4} \cline{6-8} -1 & \int\log x\,(\!\mbox{\small\it dx})^2 &
  \int\log x\,(\!\mbox{\small\it dx}) & \log x & x^{-1} &
  \mbox{-}\,x^{-2} & \mbox{\small 2!}\;x^{-3} &
  \mbox{-\small 3!}\;x^{-4} \\
  \cline{4-4} -2 & \mbox{-}\!\int\log x\,(\!\mbox{\small\it dx})
  & \mbox{-}\log x & \mbox{-}\,x^{-1} & x^{-2} &
  \mbox{-\small 2!}\;x^{-3} & \mbox{\small 3!}\;x^{-4} &
  \mbox{-\small 4!}\;x^{-5} \\
  \cline{3-3} -3 & \mbox{\small 1/2!}\; \log x & \mbox{\small 1/2!}\;
  x^{-1} & \mbox{-\small 1/2!}\;x^{-2} & x^{-3} & \mbox{-\small
    3!/2!}\;x^{-4} & \mbox{\small 4!/2!}\;x^{-5} &
  \mbox{-\small 5!/2!}\;x^{-6}\\
  \hline
\end{array}$
\smallskip
\end{center}
\caption{Tabulated results of $D^n_{\!x} x^m$}
\end{table}

\section{Riemann-Liouville (R-L) Fractional Calculus}

The R-L Fractional Calculus \cite{frac_calculus} begins with
\begin{equation}
\int_a^x f(\hat{x})\,(d\hat{x})^n\,\equiv\,
\underbrace{\int^{x}_{\!a}\!\int^{x_n}_{\!a}\!\!
  \int^{x_{n-1}}_{\!a}\!\!\!\!\!\!\!\!\cdots\!\int^{x_3}_{\!a}
  \!\!\int^{x_2}_{\!a}\!}_{n \mbox{\footnotesize-times}} f(x_1) \;
dx_1\,dx_2 \cdots dx_{n-1}\,dx_n
\end{equation}
for $n\in\Z+$ as the fundamental defining expression, and it can be
shown \cite[p. 38]{frac_calculus} to be equal to the Cauchy formula for
repeated integration,
\begin{equation}
\frac{1}{\Gamma{(n)}}
\int_a^{x}\!\!\frac{f(t)}{(x\!-\!t)^{1-n}}\,dt\;.\label{RL_int}
\end{equation}

\begin{definition}[{\bf R-L Fractional Calculus}]\quad\\
\label{d:RLFracCalc}
The R-L fractional integral is analytically extended from (\ref{RL_int}) as
\begin{eqnarray}
D^\sigma_{\!x|a} f(x)&=&\frac{d^\sigma}{dx^\sigma} f(x)
\,=\,\int_a^x f(x) (dx)^{-\sigma} \quad \mbox{by extending }
(\ref{e:Dndef})\nn\\
&=&\frac{1}{\Gamma(-\sigma)}\int_a^x\!\!
\frac{f(t)}{(x\!-\!t)^{1+\sigma}}\;dt\quad(\sigma<0,\;\sigma,a\in\R)\;\;
\mbox{by } (\ref{RL_int})\;,\label{frac_int}
\end{eqnarray}
and the R-L fractional derivative is in turn derived from the
R-L fractional integral (\ref{frac_int}) by ordinary differentiation:
\begin{equation}
D^\sigma_{\!x|a} f(x)\,=\,D^m_{\!x|a} \left( D_{\!x|a}^{-(m-\sigma)}
  f(x) \right)\quad(\sigma>0,\;m\in\Z^+)
\label{e:Dxfm}
\end{equation}
where $m$ is chosen such that $m>1+\sigma,\; \sigma>0$.
\end{definition}

\begin{lemma}
The equation (\ref{e:Dxfm}) is independent of the choice of $m$ for
$m>1\!+\!\sigma,\; m\in\Z^+,\; \sigma\in\R,\; \sigma>0$.
\end{lemma}

\noindent{\bf Proof}

For $m>1\!+\!\sigma,\; m\in\Z^+,\; \sigma>0$, we have
$-(m\!-\!\sigma)<-1<0$ and $(m-\sigma-1)>0$.  The first condition,
$-(m-\sigma)<0$, allows us to use the equation (\ref{frac_int}) to write
$$ D_{\!x|a}^{-(m-\sigma)} f(x)\,=\, \frac{1}{\Gamma(m-\sigma)}\int_a^x\!\!
\frac{f(t)}{(x\!-\!t)^{1+\sigma-m}}\;dt\;. $$
From (\ref{e:Dxfm}),
\begin{eqnarray}
\lefteqn{D^m_{\!x|a} \left( D_{\!x|a}^{-(m-\sigma)} f(x) \right)}\nn\\
&=& \frac{d^m}{dx^m} \left( \frac{1}{\Gamma(m\!-\!\sigma)}\int_a^x\!\!
\frac{f(t)}{(x\!-\!t)^{1+\sigma-m}}\;dt \right)\nn\\
&=& \frac{1}{\Gamma(m\!-\!\sigma)}\int_a^x f(t) \,
\left( \frac{d^m}{dx^m} (x\!-\!t)^{m-\sigma-1} \right) dt\;.
\label{e:lemma_choice}
\end{eqnarray}
The second condition, $(m\!-\!\sigma\!-\!1)>0$, and the condition
$m>0$ allow us to use the second case of (\ref{e:Dxfm}). Thus, 
(\ref{e:lemma_choice}) becomes
\begin{eqnarray}
\lefteqn{\frac{1}{\Gamma(m\!-\!\sigma)}\int_a^x f(t) \, 
\frac{\Gamma(m\!-\!\sigma)}{\Gamma(-\sigma)} \, (x\!-\!t)^{-(1+\sigma)}
\;dt}\nn\\
&=&\frac{1}{\Gamma(-\sigma)}\int_a^x\!\!
\frac{f(t)}{(x\!-\!t)^{1+\sigma}}\;dt\nn\\
&=&D^\sigma_{\!x|a} f(x)\nn \;.
\end{eqnarray}
\hfill\qed

When $a=0$, (\ref{frac_int}) for $f(x) = x^r$ is well-defined only for
the half plane $r>-1$. Consequently, in the R-L Fractional Calculus,
$D^\sigma_{x|a} x^r$ is well-defined only for the half plane $r>-1$.
\begin{equation}
D^\sigma_{\!x|a} x^r = \left\{\begin{array}{ll}
\ds \frac{\Gamma(1\!+\!r)}{\Gamma(1\!+\!r\!-\!\sigma)}\,
x^{r-\sigma} & (\sigma>0, \;r>-1)\\
\ds \bigg. x^r & (\sigma=0, \;\forall \; r)\\
\ds \Bigg. \frac{\Gamma(1\!+\!r)}{\Gamma(1\!+\!r\!-\!\sigma)}\,
\hat{x}^{r-\sigma}\Big|_a^x & (\sigma<0, \;r>-1)
\end{array}\right..
\label{RL_half_plane}
\end{equation}

\section{Fractional Calculus by Cauchy Integral}

The Cauchy Integral for an analytic function $f(z)$ in the complex
plane \cite[p. 120]{ahlfors} is
\begin{equation}
f^{(n)}(z_0)\,=\,\frac{\Gamma(1+n)}{2\pi i}\int_C
\frac{f(z)}{(z-z_0)^{1+n}}\,dz\;.
\end{equation}

Analytic extension of the Cauchy Integral from $n\in\Z^+$ to
$s\in\R^+$ gives an analytic extension of $D^n f(z_0)$. However, the
analytic extension is not trivial. The term $(z-z_0)^{1+\sigma}$
will become multi-valued and thus the result may depend on the choice
of branch cut and integration path.

\section{Fractional Calculus by Fourier Transform}

In the theory of Fourier Transforms,
\begin{eqnarray*}
\tilde{f}(x)&=&\int_{-\infty}^{+\infty}f(x)\,e^{ikx}\,dx\;,\\
f(x)&=&\frac{1}{2\pi}\,\int_{-\infty}^{+\infty}
\tilde{f}(x)\,e^{-ikx}\,dk\;,
\end{eqnarray*}
and
\begin{eqnarray*} 
D^\sigma_{\!x} f(x)&=&\int_{-\infty}^{+\infty}\tilde{f}(k)\,
D^\sigma_{\!x}\!\left(e^{-ikx}\right)dk\quad(\sigma\in\R)\\
&=&\int_{-\infty}^{+\infty}(-ik)^\sigma\tilde{f}(k)\;e^{-ikx}\,dk\;.
\end{eqnarray*}
This approach is often known as the pseudo-differential operator
approach. It was shown by Z\'avada \cite{frac_fourier} to be
equivalent to the Riemann-Liouville Fractional Calculus and the
Fractional Calculus by Cauchy Integral.

\section{Functional Analytic Approach}

In the functional analytic approach, an example of a functional
integral of an operator $A$ is
\begin{equation}
(-A)^a = -\,\frac{\sin a\pi}{\pi} \int_0^\infty\!
\lambda^{a-1}\,(\lambda\id - A)^{-1} A\;d\lambda \quad
(0<a<1,\;\lambda\in\R)\;.
\label{e:func_ap}
\end{equation}

$(\lambda\id - A)$ is called the kernel of the functional integral.
The evaluation of the integral with respect to real variable $\lambda$
requires various conditions on the spectrum of the operator $A$.

The analytic extension of $D$ in the functional approach is then
obtained from replacing $A$ by $D$ in (\ref{e:func_ap}).

For details of this well-developed functional analysis approach, see
\cite{op_func_anal}.

\section{Differentiating and Integrating in non-integer $s$-dimensions}

The differential of an integer $n$-dimensional function in
$n$-dimensions can be expressed as
$$
\frac{\partial}{\partial x_1}\,\frac{\partial}{\partial x_2}\cdots
\frac{\partial}{\partial x_n}\,f(x_1,x_2,\cdots,x_n)\;.
$$
The corresponding integral can be expressed as
$$
\int f(x_1,x_2,\cdots,x_n)\,d^n x
\equiv\underbrace{ \int^{x_n}\!\!\!\cdots\int^{x_2}\!\!\int^{x_1}\!\!}_{n
\mbox{\footnotesize-times}} f(x_1,x_2,\cdots,x_n)\,dx_1 dx_2\cdots dx_n\;.
$$

If $f$ is spherically symmetric, $f = f(r)$, then
\begin{equation}
\frac{\partial}{\partial x_1}\,
\frac{\partial}{\partial x_2}\cdots
\frac{\partial}{\partial x_n} f
=\frac{\partial^{n-1}}{\partial r^{n-1}}\,
\frac{\partial}{\partial\Omega}_{\!n-1}\!\!\!\!\!\! f \;\equiv
\frac{\Gamma(n/2)}{2\pi^{n/2}}\,
\frac{\partial^{n-1}}{\partial r^{n-1}} f(r)
\end{equation}
can then be analytically extended to the differential of a non-integer
$s$-dimensional function in $s$-dimensions,
\begin{equation}
\frac{\partial^{s-1}}{\partial r^{s-1}}
\frac{\partial}{\partial\Omega}_{\!s-1}\!\!\!\!\! f(r)
=\frac{\Gamma(s/2)}{2\pi^{s/2}}\,
\frac{\partial^{s-1}}{\partial r^{s-1}} f(r)
\end{equation}
where $s\in\R$ or $s\in\C$, and $\Omega$ is the $n$-dimensional solid
angle.

Similarly, the corresponding integral
\begin{eqnarray}
\int f\;d^n\!x &\equiv& \int^{\infty}_0 r^{n-1} f(r)\,dr
\int^{2\pi}_0 d\theta_1 \int^\pi_0 \sin \theta_2 d\theta_2 \,\cdots 
\int^\pi_0 \sin^{n-2}\theta_{n-1} d\theta_{n-1} \nn\\
&=& \int^{\infty}_0 r^{n-1} f(r)\,dr \int d\Omega_{n-1}
\,=\,\frac{2\pi^{n/2}}{\Gamma(n/2)} \int^{\infty}_0 r^{n-1} f(r)\,dr
\end{eqnarray}
can be analytically extended to
\begin{equation}
\int f\;d^s\!x = \frac{2\pi^{s/2}}{\Gamma(s/2)} \int^{\infty}_0
r^{s-1} f(r)\,dr\;.
\end{equation}

This method was developed by 't Hooft and Veltman \cite{hooft} in
1960's. The method was central to an important technique called
Dimensional Regularization in Quantum Field Theory where the method is
used to isolate singularities in divergent integrals.

\section{Generating Integral of the Finite Harmonic Series}

\begin{theorem}[{\bf Generating Integral of $h(n)$}]
\begin{equation}
\fbox{$\;\ds h(n) = \log x - \frac{\Gamma(1\!+\!n)}{x^n} \int_0^x
\log\hat{x}\,(d\hat{x})^n \quad (n\in\Z)\;$}
\label{e:h(n)genint}
\end{equation}
where
$$\int_a^x f(\hat{x})\,(d\hat{x})^n\,\equiv\,
\underbrace{\int^{x}_{\!a}\!\int^{x_n}_{\!a}\!\!
  \int^{x_{n-1}}_{\!a}\!\!\!\!\!\!\!\!\cdots\!\int^{x_3}_{\!a}
  \!\!\int^{x_2}_{\!a}\!}_{n \mbox{\footnotesize-times}} f(x_1) \;
dx_1\,dx_2 \cdots dx_{n-1}\,dx_n\;.$$
\end{theorem}

\noindent{\bf Proof}

We observe that $-h(n)/n!$ appears in the coefficient of the $x^n$
term when we repeatedly integrate $\log x$:
\begin{eqnarray}
\int_0^x \log\hat{x}\,(d\hat{x}) &\!\!\!=&\!\!\! x (\, \log x - 1
\,)\;,\nn\\
\int_0^x \log\hat{x}\,(d\hat{x})^2 &\!\!\!=&\!\!\! \frac{x^2}{2} (\,
\log x - \frac{3}{2} \,)\;,\nn\\
& \vdots &\nn\\
\int_0^x \log\hat{x}\,(d\hat{x})^n &\!\!\!=&\!\!\! \frac{x^n}{n!} (\,
\log x - h(n) \,)\;. \label{e:h(n)log}
\end{eqnarray}

We can prove this observation by induction:
\begin{eqnarray}
\int_0^x \log\hat{x}\,(d\hat{x})^{n+1}
&\!\!\!=&\!\!\!
\int_0^x \left(\int_0^{\hat{x}} \log\tilde{x}\,(d\tilde{x})^n \right)
(d\hat{x})\nn\\
&\!\!\!=&\!\!\!
\int_0^x \frac{\hat{x}^n}{n!} (\,\log \hat{x} - h(n) \,) (d\hat{x})\nn\\
&\!\!\!=&\!\!\!
\frac{x^{n+1}}{(n\!+\!1)!} (\,\log x - h(n) \,) - 
\int_0^x \frac{\hat{x}^n}{(n\!+\!1)!} (d\hat{x})\nn\\
&\!\!\!=&\!\!\!
\frac{x^{n+1}}{(n\!+\!1)!} \left( \log x - h(n) - \frac{1}{n\!+\!1} 
\right)\nn\\
&\!\!\!=&\!\!\!
\frac{x^{n+1}}{(n\!+\!1)!} (\,\log x - h(n\!+\!1) \,)\;.\nn
\end{eqnarray}

Rearrangement of (\ref{e:h(n)log}) yields the Theorem.
\hfill\qed

\begin{theorem}[{\bf Generating Integral of $h(\rho)$}]
\begin{equation}
\fbox{$\;\ds h(\rho) = \log x - \frac{\Gamma(1\!+\!\rho)}{x^\rho}
\int_0^x \log\hat{x}\,(d\hat{x})^\rho \quad (\rho\in\R)\;$}\;.
\label{e:h(rho)intlog}
\end{equation}
\end{theorem}

\noindent{\bf Proof}

By analogy to generating functions, we take
$$\int_0^x \log \hat{x}\,(d\hat{x})^n$$
as the {\em generating integral} of the finite harmonic series $h(n)$,
and so the natural analytic extension of the generating integral
takes the form of
\begin{equation}
\int_0^x \log\hat{x}\,(d\hat{x})^\rho =
\frac{x^\rho}{\Gamma(1\!+\!\rho)} \left( \log x - h(\rho) \right)\;.
\label{e:intlogh(rho)}
\end{equation}

Noting that $\log x$ may be expressed as
\begin{equation}
\log x = \int^x_1 \hat{x}^{-1} d\hat{x} = \lim_{\epsilon \to 0} \int_1^x
\hat{x}^{-1+\epsilon} \, d\hat{x} = \lim_{\epsilon \to 0}
\frac{1}{\epsilon}\,( x^\epsilon - 1 )\;,
\end{equation}
we can now evaluate the fractional integral in (\ref{e:intlogh(rho)}) by the
{\bf R-L Fractional Calculus} (Definition \ref{d:RLFracCalc}).
\begin{eqnarray} \int_0^x \log\hat{x}\,(d\hat{x})^\rho 
&\!\!\!=&\!\!\! \int_0^x \lim_{\epsilon
  \to 0} \frac{1}{\epsilon}\,(\hat{x}^\epsilon - 1)\,
  (d\hat{x})^\rho
\,=\, \lim_{\epsilon \to 0} \frac{1}{\epsilon} \int_0^x
  (\hat{x}^\epsilon - 1)\, (d\hat{x})^\rho \nn\\
&\!\!\!=&\!\!\! \lim_{\epsilon \to 0} \frac{1}{\epsilon}
  \left[ \int_0^x \hat{x}^\epsilon (d\hat{x})^\rho - \int_0^x 1\,
  (d\hat{x})^\rho \right]\nn\\
&\!\!\!=&\!\!\! \lim_{\epsilon \to 0} \frac{1}{\epsilon}\,
  \bigg[ D^{-\rho}_{\!\hat{x}} \hat{x}^\epsilon \Big|_0^x \;\,-\;\,
  D^{-\rho}_{\!\hat{x}} \,1\, \Big|_0^x \bigg] \nn\\ 
&\!\!\!=&\!\!\! \lim_{\epsilon \to 0} \frac{x^\rho}{\epsilon} \left[ 
\frac{\Gamma(1\!+\!\epsilon)\,x^\epsilon}
{\Gamma(1\!+\!\epsilon\!+\!\rho)} -
\frac{1}{\Gamma(1\!+\!\rho)} \right]\label{e:loggamma}\end{eqnarray}
where the interchange of the integral and the limit is justified by
Arzel\`a's theorem on bounded convergence \cite[pp. 405-406]{arzela} as 
$(\hat{x}^\epsilon-1)/\epsilon$ is integrable in $\hat{x}\in
[0,x]\,$.

Combining (\ref{e:loggamma}) with (\ref{e:intlogh(rho)}) gives the
analytic extension of the finite harmonic series:
\begin{eqnarray} h(\rho)
&\!\!\!=&\!\!\! \ds \log x - \lim_{\epsilon \to 0} \frac{1}{\epsilon}
\left[ \frac{\Gamma(1\!+\!\epsilon)\,
    \Gamma(1\!+\!\rho)}{\Gamma(1\!+\!\epsilon\!+\!\rho)}\,
x^\epsilon - 1\right]\nn\\
&\!\!\!=&\!\!\! \ds \log x - \lim_{\epsilon \to 0}\,
\frac{1}{\epsilon}\,(x^\epsilon - 1) \left[ \frac{\Gamma(1\!+\!\epsilon)\,
    \Gamma(1\!+\!\rho)}{\Gamma(1\!+\!\epsilon\!+\!\rho)} \right]\nn\\
&&\ds\quad
+ \;\lim_{\epsilon\to 0} \frac{1}{\epsilon} \left[ 1 -
  \frac{\Gamma(1\!+\!\epsilon)\,
    \Gamma(1\!+\!\rho)}{\Gamma(1\!+\!\epsilon\!+\!\rho)} \right]\nn\\
&\!\!\!=&\!\!\! \ds \left[\log x - \lim_{\epsilon \to 0}\,
\frac{1}{\epsilon}\,(x^\epsilon - 1)\right]
+\lim_{\epsilon\to 0}\frac{1}{\epsilon} \left[ 1 -
  \frac{\Gamma(1\!+\!\epsilon)\,
    \Gamma(1\!+\!\rho)}{\Gamma(1\!+\!\epsilon\!+\!\rho)} \right]\nn\\
&\!\!\!=&\!\!\! \ds \lim_{\epsilon \to 0} \frac{1}{\epsilon} \left[ 1
  - \frac{\Gamma(1\!+\!\epsilon)\,
    \Gamma(1\!+\!\rho)}{\Gamma(1\!+\!\epsilon\!+\!\rho)} \right]\nn\\
\Bigg.&\!\!\!=&\!\!\!\ds\frac{\Gamma'(1\!+\!\rho)}{\Gamma(1\!+\!\rho)} -
\Gamma'(1)\nn\\
\Big.&\!\!\!=&\!\!\!\ds\psi(1\!+\!\rho) + \gamma\label{e:harmonic}
\end{eqnarray}
where the limit has been taken with L'Hospital rule.
\hfill\qed

\begin{figure}[hbt]
\begin{center}
\begin{tabular}{c}
\raisebox{50pt}{$\!\!\!\!\!\!h(\rho)$}
\includegraphics[height=125pt]{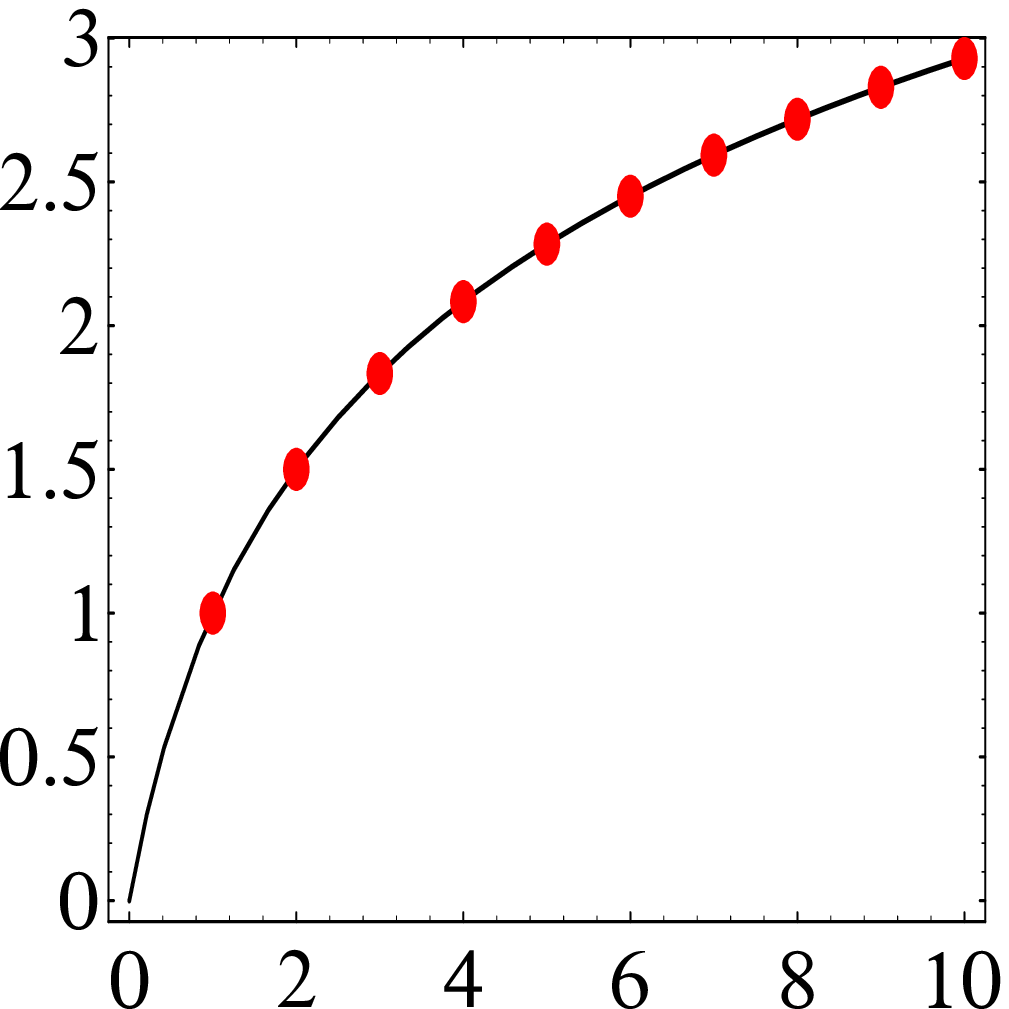}\\
$\qquad\rho$
\end{tabular}
\caption{The curve $\ds h(\rho)=\psi(1+\rho)\!+\!\gamma$
  passes through the points $\ds (n,h(n))$ where $n\in\Z^+, \;
  h(n)=\ds \sum_{k=1}^n\frac{1}{k}\;$.}
\end{center}
\end{figure}

Hence, by the application of the R-L Fractional Calculus to
analytically extend the generating integral, we have found an
alternative elementary proof of the first part (\ref{e:harmthm1}) of
the Finite Harmonic Series Theorem.

\begin{theorem}
\label{t:Dlog}
\begin{equation}
\fbox{$\;\ds \int_0^x \log\hat{x}\,(d\hat{x})^\rho =
\frac{x^\rho}{\Gamma(1\!+\!\rho)} \left( \log x -  \psi(1\!+\!\rho) -
  \gamma \right) \quad (\rho\in\R)\;$}\;.
\label{e:Dlog}
\end{equation}
\end{theorem}

\noindent{\bf Proof}

Replacing the $h(\rho)$ in (\ref{e:intlogh(rho)}) by (\ref{e:harmonic})
gives the Theorem.
\hfill\qed

\section{The Riemann Zeta Function up to Order $n$}

The analytic extension (\ref{e:harmonic}) can be generalized.

\begin{definition}
The Riemann zeta function \cite{titchmarsh}
\begin{equation}
\zeta(s)\,=\,\sum_{k=1}^{\infty}\frac{1}{k^s}\quad (\re(s)>1),
\end{equation}
the polygamma functions \cite[p. 260, (6.4.1)]{abramowitz}
\begin{equation}\psi^{(m)}(x)\;=\;\frac{d^m}{dx^m}\psi(x)\;=\;
\frac{d^{m+1}}{dx^{m+1}} \log \Gamma(x)
\,=\,(-1)^{m+1} m! \,\sum_{k=0}^{\infty}
\frac{1}{{(x\!+\!k)}^{m+1}}\;,\end{equation}
and the Riemann zeta function up to order $n$
\begin{equation}
\zeta(s \bigm| n)\;=\;\sum_{k=1}^n\frac{1}{k^s}\;=\;\zeta(s)
-\!\!\sum_{k=n+1}^{\infty}\frac{1}{k^s}\quad (\re(s)>1)
\end{equation}
may be combined to write
\begin{eqnarray} \zeta(m \bigm| n)&=&\sum_{k=1}^n\frac{1}{k^m}
\,=\,\frac{(-1)^m}{(m\!-\!1)!}\left(\psi^{(m-1)}(1\!+\!n)-\psi^{(m-1)}(1)
\right)\nn\\
&=&\frac{(-1)^m}{(m\!-\!1)!}\frac{d^m}{dx^m}\,\log(\Gamma(1\!+\!x))\,
\Big|_{x=0}^{x=n}\;.\label{e:zetadlog}\end{eqnarray}
\end{definition}

The analytic extension is then obtained by replacing
the derivative in (\ref{e:zetadlog}) with a fractional
derivative:
\begin{equation}\zeta(s \bigm| z)\;\;=\;\frac{w(s)}{\Gamma(s)}\;
D^s_{\!x}\,\log(\Gamma(1\!+\!x))\,\Big|_{x=0}^{x=z}\quad (s,z\in\C)
\label{e:zetadloggamma}
\end{equation}
which can be evaluated when $\log(\Gamma(1\!+\!x))$ is expressed in
the form of an asymptotic series \cite[p. 257, (6.1.41)]{abramowitz}.
However, $w(s)$ depends on the choice of extension to the {\bf R-L
  Fractional Calculus} (Definition \ref{d:RLFracCalc}) into the other
half plane, $r\le-1,\;\sigma\in\R\,$.

\section{Analytic Extension of R-L Fractional Calculus}

Consider the case of $D^n_{\!x}$ where $a = 0$.  We shall introduce
the $(\sigma,r)$ diagram in which the numerical factor of
$D^\sigma_{\!x} x^r$ is mapped to the point at coordinate $(\sigma,r)$
of the diagram.  The $(\sigma,r)$ diagram of $D^\sigma_{\!x} x^r$ can
be characterised into 4 regions as in Figure \ref{fig:sigma_r_diagram}:

\begin{figure}[hbt]
\begin{center}
\includegraphics[height=180pt]{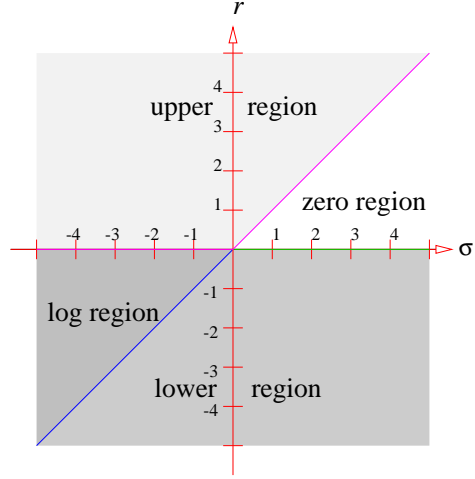}
\caption{$(\sigma,r)$ diagram of $D^\sigma_x x^r.$}
\label{fig:sigma_r_diagram}
\end{center}
\end{figure}

\begin{definition}[{\bf Regions of $(\sigma,r)$}]\quad\\
\label{d:regions}
\begin{tabular}{lll}
The zero  region&$Z_{ero}$&$= \{(\sigma,r) : \;r<  \sigma, \;r\ge0\}$;\\
the upper region&$U_{pp} $&$= \{(\sigma,r) : \;r\ge\sigma, \;r\ge0\}$;\\
the lower region&$L_{ow} $&$= \{(\sigma,r) : \;r<  \sigma, \;r<  0\}$;\\
the log   region&$L_{og} $&$= \{(\sigma,r) : \;r\ge\sigma, \;r<  0\}$.
\end{tabular}
\end{definition}

A point lying on the right of the $r$-axis $(\sigma>0)$ is a
differentiation; a point on the left $(\sigma<0)$ is an integration.

The $(\sigma,r)$ diagram at integer grid points everywhere except in
the log region gives numerical factors identical to those in Table 1.

An extension of R-L Fractional Calculus to the other half plane $r\le
-1$ is given by
\begin{equation}
D^\sigma_{\!x} x^r = \left\{
\begin{array}{ll}
\ds \lim_{\epsilon\to 0}\,
\frac{\Gamma(1\!+\!r\!+\!\epsilon)}{\Gamma(\epsilon)}\,
D^{\sigma-r}_{\!x} \log x & \mbox{ in } \Omega\\
\ds \lim_{\epsilon\to 0}\,\frac{\Gamma(1\!+\!r\!+\!\epsilon)}
{\Gamma(1\!+\!r\!+\!\epsilon\!-\!\sigma)}\,x^{r-\sigma}
& \mbox{ elsewhere}
\end{array}\right.\label{D_combined}
\end{equation}
where $\Omega = \{(\sigma,r) : \;\sigma\in\R, \;r\in\Z^-\}$, the set
of horizontal lines in lower and log regions.

$\Gamma(1\!+\!r)/\Gamma(1\!+\!r\!-\!\sigma)$ is finite everywhere in
the zero and upper regions.  $\stack{{}_{\mbox{\small
      lim}}}{\epsilon\to 0} \;
\Gamma(1\!+\!r\!+\!\epsilon)/\Gamma(1\!+\!r\!+\!\epsilon\!-\!\sigma)$
is well-defined everywhere in the lower and log regions except in 
$\Omega\backslash(\Z^-\times\Z^-)$. Following Theorem \ref{t:Dlog},
the R-L fractional integral of $\log x$ can be evaluated exactly and
expressed in only elementary functions,
\begin{equation}
D^\sigma_{\!x}\log x = \left\{ \begin{array}{ll}
\ds \frac{x^{-\sigma}}{\Gamma(1\!-\!\sigma)} \left(\big.\log x - 
\psi(1\!-\!\sigma)-\gamma\right) 
& (\sigma\in\R\big\backslash\Z)\\
\ds\Bigg. \lim_{\epsilon\to 0}\, \frac{\Gamma(1\!+\!r\!+\!\epsilon)}
{\Gamma(1\!+\!r\!+\!\epsilon\!-\!\sigma)}\,x^{r-\sigma} & (\sigma\in\Z)
\end{array} \right..
\label{log_psi}
\end{equation}
(\ref{D_combined}) is thus well-defined.

To analytically extend from $D^\sigma_{\!x} x^r$ on the real plane
$(\sigma,r)$ to $D^s_{\!x} x^r$ on the product of complex plane and
real line, $(s,r)\in\C\times\R$, we simply replace $\sigma\in\R$ in
(\ref{D_combined}) by $s\in\C$.

\section{Open Problems}

\begin{enumerate}

\item Generalize (\ref{D_combined}) for $D^s_{\!z-c} z^w, \;s,w,z,c\in\C$.
For the complex function $z^w$, one has to specify, in addition, the
integration contour for $\re(s)<0$.

\item Find the exact expression for $\log(\Gamma(1\!+\!x))$ and $w(s)$ in
terms of elementary functions in analytic extension of the R-L
Fractional Calculus given by (\ref{D_combined}).

\end{enumerate}

\section{Tables of Generating Integrals}

\begin{quote}{\em ``Nature laughs at the difficulties of
    integration.''} --- Laplace\end{quote}

\smallskip

An interesting consequence is that objects of the form $\int x^r (\log
x)^a (dx)^\rho$, $r,a,\rho\in\R$, exist and can be generating integrals
for certain functions. Perhaps it may be worthwhile to introduce, in
the future editions of Tables of Integrals, a new section which gives
the coefficient of the $x^k$ term, $w(\rho,r,a,k)$, corresponding to
these generating integrals to facilitate the evaluation of integrals
of similar forms.

\end{document}